\documentclass[12pt]{article}
\usepackage{theorem,amsmath,array,amsfonts,amssymb}

\newenvironment{proof}{{\bf Proof}:\quad
                     }{{\hfill$\Box$\\}}

\begin{document}
\begin{center}
{\bf \LARGE Some Remarks on the Local Fundamental Group Scheme}
\end{center}
\vskip 1.5cm
\begin{center}
{\bf V.B. Mehta and S. Subramanian} \\
School of Mathematics\\ Tata Institute of
Fundamental Research \\ Homi Bhabha Road \\ Mumbai 400 005 \\ India. \\
[5mm]
e-mail: vikram@math.tifr.res.in \\
~~~~~~~~~~~subramnn@math.tifr.res.in
\end{center}

\vskip 5mm
\begin{center}
{\Large \bf Introduction}
\end{center}
\vskip 1cm

The fundamental group - scheme was introduced by M. Nori in [3,4]. In the
same paper he had made 2 conjectures [4]. The first conjecture was proved
 by the present authors in [2].  In  an attempt to prove the second
conjecture,
 we had introduced the``local fundamental group-scheme'' which is defined only 
in characteristic $p$.  In this note, we gather together some necessary and 
sufficient conditions for it to base-change ``correctly'' from one
algebraically closed field to another algebraically closed field of
characteristic $p$.

\section*{Section ~~ I}

We define the local fundamental group scheme, denoted by $\pi^{loc}$,as
follows:

   Let $X$ be a smooth projective variety over an algebraically closed field of characteristic $p$, and let $F: X \to X$ be the Frobenius map.  Let
$FT(X)$ denote the category of all vector bundles $V$ on $X$ such that
$F^{t *} (V)$ is {\it trivial} for some integer $t$.  Fix a base point
$x \in X$ and let $T:FT(X)$  to Vect $k$ be given by $V \to 
V_x$.  One checks easily that $FT(X)$, with the fibre functor $T$, is a 
Tannaka category [1].  The corresponding affine group - scheme is defined
to 
be the {\it local fundamental group-scheme} of $X$, denoted by
$\pi^{loc}(X)$. Note that $\pi^{loc}(X)$  is defined only in characteristic $p$.

\paragraph*{Remark 1}  We have a canonical surjection $\pi (X) \to   
 \pi^{loc}(X)$, where $\pi (X)$ is the fundamental group scheme as
defined by Nori in [4].

\paragraph*{Remark 2}  Recall that all varieties are smooth and projective, and defined over algebraically closed fields of characteristic $p$.  Now we have 

\paragraph*{Lemma 1.1}a)  If $V \in F T(X)$, then $V$ is semistable for 
any 
polarization $H$ on $X$.

b)  If $V$ is stable for 1 polarization, then $V$ is stable with respect to any polarization.

\paragraph*{Proof a} If $F^{t*}(V)$ is a trivial vector bundle on $X$
for
some $t$, then it is semistable for any polarization.  It follows that 
$V$ itself is semistable with respect to any polarization.

b)  If  $V \in FT (X)$, then there exist
\begin{enumerate}
\item[1)] a principal $G$-bundle \\
$$E \to X$$
where  $G$ is a finite local group-scheme,and E is {\it reduced} [4]
\item[2)]
\end{enumerate}

and a representation $$\sigma:G \to GL (n), n = \ {\rm rank} \ V,$$ 
such that $V \cong V_\sigma$.  It follows easily that subbundles $W$ of $V$
with $\mu (W) = \mu (V)$ are given by $G$ subspaces $W_1$ of $k^n$, hence 
if $V \in FT(X)$ and stable with respect to one polarization, then it is 
stable with respect to any polarization.

\paragraph*{Definition 1.2:}   Fix positive integers $r$ and $t$.  Let
$S (X, r,t)$ denote the {\it set of isomorphism classes} of stable vector 
bundles $V$ on $X$ such that rank $V=r$ and $F^{t*} (V)$  is trivial.
We shall show that $S(X, r,t)$ corresponds bijectively to the  {\it closed}
points of a scheme $M(X, r,t)$ defined over $k$. 

\paragraph*{Remark 3}:We recall that all fields  are algebraically closed
of characteristic $p$ and all curves are assumed smooth and projective.

\paragraph*{Theorem}:The following statements are equivalent.

1) $$g: \pi^{\rm loc} (X_{k'}, x) \to \pi^{\rm loc} (X_k, x)\times_k \ 
{\rm 
spec} k'.$$ is an isomorphism , for every $k' \supset k$

2) $S(X;r,t)$ is a finite set , for all $r$ and $t$.

3)For every $k' \supset k$, every $V \in S(X_{k'};r,t,)$ descends
to $X$, i.e. is defined over $k$.

\begin {proof}
Before we prove that the statements are equivalent, we discuss in brief
some generalities on Tannaka categories.
Let ${\cal C}'$ be the Tannaka category of $F$-trivial vector bundles on 
$X_k$.
Then ${\cal C}'$ is isomorphic to the Tannaka category of rational, finite dimensional representations of $\pi^{loc} (X,x)$.  We consider the category 
${\cal C}_1'$ defined as follows:

a vector bundle $V$ on $X_{k'}$ is an object of ${\cal C}_1'$ if there is
an $F$-trivial vector bundle $V_1$ on $X_k$ such that $V$ is a degree zero 
subquotient of $V_1 \otimes_k k'$ on $X_{k'}$.  The morphisms of ${\cal 
C}_1'$ are vector bundle homomorphisms and it is easily seen that with the tensor
product of vector bundles as the tensor operation in the category and the
fibre functor which associates to a bundle $V$ its fibre over the base point $x$ of $X_{k'}, {\cal C}_1'$ is a Tannaka category.

Let $R'$ denote the Tannaka category of rational, finite dimensional
representations of $\pi^{\rm loc}(X_k, x)\times_k  spec k'$ over $k'$.  
We can 
define a functor from $R'$ to ${\cal C}_1'$ as follows:

If $V$ is an object of $R'$ namely a representation of $\pi^{\rm loc}(X_k, 
x) \times_k  Spec k'$,  then there is a finite, local group scheme $G$ 
over $k$ and a surjection $\pi^{\rm loc} (X_k, x) \to G$ such that $V$ is actually a 
representation of $G_{k'} =G \times_k $spec $k'$.  The homomorphism 
$\pi^{\rm loc} \ (X_k,
x) \to G$ defines a principal $G$-bundle $\pi: E \to X_k$ on $X_k$.  Let
$\pi: E_{k'} \to X_{k'}$ be the base change of this bundle to $k'$ so that
$E_{k'} \to X_{k'}$ is a principal $G_{k'}$ bundle.  The representation
$V$ of $G_{k'}$ associates to this principal bundle a vector bundle (which 
is $F$-trivial) denoted by $\tilde{V}$.  The $G_{k'}$ module $V$ is a 
subquotient
of $k'[G]^{\oplus n}$ where  $k' [G]$ denotes the coordinate ring of 
$G_{k'}$ (see).  As $k' [G] = k [G] \otimes_k k'$ where $k[G]$ denotes the coordinate ring of $G$ over $k$, if we let $V_1$ be the vector bundle on $X_k$ associated
 to the $G$-bundle $\pi:E \to X_k$ by the representation of $G$ on 
$k[G]^{\oplus n}$, we obtain that $\tilde{V}$ is a degree zero subquotient of
$V_1 \otimes_k k'$.  Hence $\tilde{V}$ is an object of ${\cal C}_{1}'$.  It is 
easily seen that the functor from $R'$ to ${\cal C}_1'$ defined as above 
is 
an equivalence of Tannaka categories, so the affine group scheme associated
to ${\cal C}_1'$ is $\pi^{\rm loc} \ (X_k,x) \times_k$ spec $k'$.

Let ${\cal D}'$ denote the Tannaka category of all $F$-trivial vector bundles 
on $X_{k'}$.  Since every object of ${\cal C}'_1$ is $F-$trivial ${\cal C}_1'$
is a Tannaka subcategory of ${\cal D}'$.  The affine group scheme 
associated
to ${\cal D}'$ is $\pi^{\rm loc} \ (X_{k'},x)$ and the functor ${\cal
C}_1'
C {\cal D}'$ induces the canonical homomorphism
$$g: \pi^{\rm loc} (X_{k'}, x) \to \pi^{\rm loc} (X_k, x) \times_k \ {\rm 
spec} k'.$$
An application of Proposition 5, p.121 in [4] shows that $g$ is 
surjective. Now we prove 1) ,2) and 3) are equivalent.

\paragraph*{$1 \Rightarrow 3$} If $g$ is an isomorphism, then it is, in
particular, a closed immersion, and arguing as in Proposition (3.1) of
[2],
we see that any stable $F$-trivial bundle over $X_{k'}$ descends to $X_k$.

\paragraph*{$3 \Rightarrow 1$}  We now assume that every stable $F$-trivial bundle on $X_{k'}$ is defined over $X_k$.  We have to show that $g$ is a closed immersion.  We use Proposition (2.21) (b) in [] for this purpose.  We have to show that any object

of ${\cal D}'$, namely an $F$-trivial bundle on $X_{k'}$, is a subquotient
of an object of ${\cal C}_1'$.  Since any $F$-trivial vector
bundle on $X_{k'}$ is semistable, it has a filtration by stable $F$-trivial bundles on $X_{k'}$.  By  hypothesis, the stable $F$-trivial bundles descend to $X_k$, and we are left  with an extension of stable bundles.  For simplicity, let
$V$ be an object of ${\cal D}'$, which is an extension.
$$0 \to V_1 \to V \to V_2 \to 0$$
where $V_1 $ and $V_2$ are stable $F$-trivial bundles on $X_{k'}$.  The
extension bundle $V$ defines an element of $Ext^1_{{\cal O}_{X_{k'}}} (V_2,V_1)
$.  However, we have

$$Ext^1_{{\cal O}_{X_{k'}}} (V_2, V_1) = Ext^1_{{\cal O}_{X_{k}}} (V_2, V_1) \otimes_k k'.$$
Thus, the element representing $V$ in $Ext^1_{{\cal O}_{X_{k'}}} (V_2, V_1)$ 
can be written as $\sum_i \alpha_i e_i$, where $\alpha_i \in k'$ and
$e_i \in Ext^1_{{\cal O}_{X_k}}(V_2, V_1)$.  By the construction of Baer
sum in $Ext^1$ (see Definition (3.4.4), p.78, in [5]), we see that $\sum \alpha_i e_i$ represent a subquotient of extensions of $V_2$ by $V_1$ on $X_k$.  Hence
we obtain that $V$ is a subquotient of a $F$-trivial bundle defined over $X_k$.
 The case where there are more than two stable factors in the stable
filtration of $V$ is similar.
\paragraph*{2 $\Leftrightarrow$ 3}  Let $O_t$ be the open subset of $R^s$ 
(for the definition of $R^s$, see $\S $4, p.141 in [6]) defined by \\
$O_t = \{x \in R^s \mid F^{t*} (U_x)$ is  semistable  on $X$ where $F$ is
the Frobenius of $X$\} \\ .  Let $Y_t$ be the closed subset of $O_t$
defined by \{$y \in O_t \mid F^{t*} U_y$ is trivial\} 
and we give it the reduced subscheme structure.
As $O_t$ is open and $GL(N)$ invariant in $R^s$, the quotient $O_t /GL (W)$
exist as a scheme.  Since $Y_t$ is closed and $GL(N)$ invariant in $O_t$, the
quotient $Y_t / GL(N)$ exists as a scheme, and let this scheme be denoted by
$M(X,r,t)$.  The $k$-rational points of $M(X,r,t)$ are in bijection with 
$S(X,r,t)$.  It is  further clear that $k'\supset k$ the schemes $O_t, 
Y_t$ and
$M(X,r,t)$ base change to the corresponding schemes over $k'$.  In particular,
$$M(X_{k'}, r,t) \cong M (X_k, r,t) \times_k {\rm spec} \ k'.$$
In general, one always has
$$M(X_{k'},r,t) (k) \subset M (X_{k'},r,t) (k').$$
If strict equality holds, 
$M(X,r,t)$ is the spectrum of an Artin ring on $k$.
In that case, every element of $S(X_{k'}, r,t)$ is  defined over $k$.  Further,
we see that if every element of $S(X_{k'}, r,t)$ is defined over $k$, then $S
(X,r,t)$ is finite.
\end{proof}

\begin{center}
{\bf References}
\end{center}
\begin{enumerate}
\item[{[1]}] P. Deligne, J. Milne:  Tannakian categories in Hodge cycles, 
Motives and Shimura varieties, {\it Lecture Notes in Mathematics} 900, 
Springer.
\item[{[2]}] V.B. Mehta, S. Subramanian: On the fundamental group scheme, 
{\it Inv. Math.} 148, 143-150 (2002).
\item[{[3]}] M.V. Nori, Representations of the fundamental group, 
{\it Compositio Math.} 33, 1976.
\item[{[4]}] M.V. Nori, The fundamental group scheme, {\it Proc. Indian Acad. 
Sci.} (Math. Sci) 91, 1982.
\item[{[5]}] Charles A. Weibel, An introduction to homological algebra, 
{\it Cambridge university Press}, 1994.
\item[{[6]}] P.E. Newstead, Introduction to Moduli problems and Orbit 
spaces,
{\it TIFR Lecture Notes}, 1978.
\end{enumerate}

\end{document}